\documentclass[11pt]{amsart}

\usepackage{hyperref}

\usepackage{bbm,stmaryrd,mathrsfs,mathabx,amssymb}
\usepackage[all]{xy}

     \def\Si{\Sigma}
       \def\th{\theta}
\def\al{\alpha}       
\def\be{\beta}        
\def\ga{\gamma}       
      \def\ka{\kappa}
\def\vi{\varphi}

\def\gA{\mathfrak a} \def\gN{\mathfrak n}
\def\gB{\mathfrak b} 
 \def\gP{\mathfrak p}
 
 \def\gR{\mathfrak r}

\def\gK{\mathfrak k}

 \def\dP{\mathfrak P}

\def\dL{\mathfrak L}

\def\kA{\mathcal A} 
 
\def\kC{\mathcal C} \def\kP{\mathcal P}

 \def\kX{\mathcal X}
 
\def\kM{\mathcal M}

 \def\cX{\mathscr X}

\def\tMd{\mbox{-}\widetilde{\mathrm{Mod}}}
\def\tM{\tilde{M}}		\def\ti{\tilde}
\def\ox{\bar{x}}		\def\hA{\hat{A}}

\def\bap{\bigcap}        \def\bup{\bigcup}

\def\*{\otimes}		\def\+{\oplus}		
\def\bop{\bigoplus}	
\def\sb{\subset}         
\def\spe{\supseteq}      \def\sbe{\subseteq}
\def\0{\emptyset}		\def\8{\infty}

\def\xarr{\xrightarrow}
\def\emb{\hookrightarrow}

\def\<{\langle}		\def\>{\rangle}

\def\Arr{\Rightarrow}

\def\={{\smallsetminus}}

\def\smtr#1{\left(\begin{smallmatrix}#1\end{smallmatrix}\right)}

\def\lst#1#2{ #1_1 , #1_2 , \dots , #1_{#2} }
\def\wrd#1#2{ #1_1  #1_2  \dots  #1_{#2} }

\def\Hom{\mathop\mathrm{Hom}\nolimits}
\def\thom{\mathop{\widetilde{\mathrm{Hom}}}\nolimits}

\def\Ext{\mathop\mathrm{Ext}\nolimits}
\def\tExt{\mathop{\widetilde{\mathrm{Ext}}}\nolimits}

\def\gdim{\mathop\mathrm{l.gl.dim}\nolimits}

\def\ker{\mathop\mathrm{Ker}\nolimits}

\def\im{\mathop\mathrm{Im}\nolimits}

\def\ann{\mathop\mathrm{ann}\nolimits}
\def\rad{\mathop\mathrm{rad}\nolimits}

\def\cok{\mathop\mathrm{Coker}\nolimits}
\def\he{\mathop\mathrm{ht}}
\def\trs{\mathop\mathrm{tors}}
\def\tf{\mathop\mathrm{tf}}
\def\soc{\mathop\mathrm{soc}\nolimits}
\def\idim{\mathop\mathrm{inj.dim}\nolimits}
\def\gdim{\mathop\mathrm{gl.dim}\nolimits}

\def\id{\mathrm{Id}}

\def\Md{\mbox{-}\mathrm{Mod}}

\def\ms{S^{-1}}

\def\pn{pseudo-noetherian}

\def\pno{pseudo-noetherian order}
\def\pk{pseudo-krullian}
\def\pko{pseudo-krullian order}
\def\qd{quasidivisorial}
\def\nzd{non-zero-divisor}

\newtheorem{thm}{Theorem}[section]
\newtheorem{lem}[thm]{Lemma}
\newtheorem{prop}[thm]{Proposition}
\newtheorem{cor}[thm]{Corollary}
\newtheorem{claim}{Claim}[thm]
\newtheorem{fact}[thm]{Fact}

\theoremstyle{definition}
\newtheorem{dfn}[thm]{Definition}
\newtheorem{exam}[thm]{Example}

\theoremstyle{remark}
\newtheorem{remk}[thm]{Remark}

\numberwithin{equation}{section}

\begin{document}

 \title{Injective modules over pseudo-krullian orders}
 \author{Yuriy A. Drozd}
 \address{Harvard University and Institute of Mathematics of the National Academy of Sciences of Ukraine}
 \email{y.a.drozd@gmail.com}
 
 
  \begin{abstract}
  We introduce a new class of rings, \emph{pseudo-krullian orders}, consider the Serre quotients of their module categories with respect to pseudo-isomorphisms and describe injective objects in such quotient categories and its global homological dimension. These results generalize the results of I. Beck for the case of Krull rings. In particular, we establish the global homological dimension of the category of maximal Cohen-Macaulay modules over an order over a noetherian ring of Krull dimension 2.
  \end{abstract}
  
 \maketitle
 
  Theory of divisors, originated from the classical papers of Kronecker on algebraic numbers, developed into a vast part
  of commutative algebra and algebraic geometry. The most accomplished form it has for \emph{Krull rings}, in particular,
  normal (integrally closed) noetherian rings (see, for instance, \cite[Ch.\,VII]{bourbaki}). In particular, in this case the 
  reduction of ideals to divisors is naturally extended to arbitrary modules as \emph{pseudo-isomorphism},
  that is ``isomorphism in codimension $1$'' (ibid.). Actually, it is a special case of Serre quotient for the category of modules,
  which clarify some questions about structure of modules. In his paper \cite{beck} I.\,Beck has studied the structure
  of this quotient with special attention to injective objects and their relations with injective modules. In this paper we 
  generalize his results to much more general situation of \emph{\pko{s}}, maybe noncommutative and, if commutative domains,
  not necessarily integrally closed. Following I.\,Beck, we introduce a special class of modules (\emph{codivisorial modules}) and
  establish its relations with the Serre quotient with respect to pseudo-isomorphisms. We also show that injective codivisorial
  modules behave just as injective modules over noetherian rings and describe injective objects in the 
  quotient category. As a corollary, we establish the global homological dimension of this quotient category. 
  As an application, we find the global homological dimension of  the category of maximal Cohen-Macaulay 
  modules over a finite algebra over a noetherian ring of Krull dimension $2$. In particular, in case of a normal domain or some
  orders over such domain this category is hereditary (of global homological dimension $1$).
  
  \medskip  Some remarks about notations.
  We denote by $\sb$ the \emph{proper embedding} (``less''), so $\sbe$ means ``less or equal''.
  We write ``iff'' instead ``if and only if''. All rings are supposed \emph{unital}, that is having a unit, all ring homomorphisms 
  mapping unit to unit and all modules $M$ \emph{unital}, that is such that $1x=x$ for every $x\in M$.

 \section{Pseudo-krullian orders and codivisorial modules}
 \label{sec1} 

 Recall that a commutative ring $R$ is called \emph{\pn} \cite{sem,cher} if the following conditions hold:
 \begin{enumerate}
  \item For every element $a\in R$ the set $V_{min}(a)$ of prime ideals minimal among those containing $a$ is finite.
 \item  For every $\gP\in V_{min}(a)$ the ring $R_\gP$ is noetherian.
  \end{enumerate} 
  Noetherian rings and Krull rings are examples of \pn\ rings. One easily sees that $\he\gP\le1$ for every prime ideal 
  $\gP\in V_{min}(a)$ and $\he\gP=1$ if $a$ is a \nzd. Let $\{\lst\gN s\}$ be the set of minimal prime ideals of $R$
  (that is $V_{min}(0)$). Then the rings $R_{\gN_j}$ are artinian. If, moreover, $R$ is reduced (has no nilpotent elements), 
  $Q_j=R_{\gN_j}$ is a field and $R$ embeds into the semisimple ring $Q=\prod_{j=1}^sQ_j$, which is the full ring of
  fractions of $R$. In this case $Z(R)=\bup_{j=1}^s\gN_j$ is the set of zero divisors of $R$ and $\bap_{j=1}^s\gN_j=0$.
  For an $R$-module $M$ we write $QM$ instead of $Q\*_RM$. We say that an element $x\in M$ is \emph{torsion}
  if there is a \nzd\ $a\in R$ such that $ax=0$ and \emph{torsion free} otherwise. Obviously, $x$ is torsion iff $1\*x=0$
  in $QM$. If all elements of $M$ are torsion (resp., torsion free), we call the module $M$ \emph{torsion} (resp., torsion
  free), and we identify a torsion free module with its image in $QM$. We denote by $\trs M$ the submodule of all torsion
  elements and by $\tf M$ the quotient $M/\trs M$. 
  
  An algebra $A$ over a \pn\ ring $R$ is called \emph{\pn} if $A_\gP$ is noetherian for every prime ideal $\gP$ with $\he\gP\le1$.
  If $R$ is reduced, $QA$ is a finite algebra over the artinian ring $Q$; if $A$ is reduced itself (has no nilpotent ideals) this algebra 
  is semisimple. If, moreover, $A$ is torsion free as $R$-module, hence embeds into $QA$, we call it a \emph{\pno}.
  
  \vskip.5ex
  In what follows $R$ is a reduced \pn\ ring and $A$ is a reduced \pno\ over $R$. We denote by $\kP=\kP(R)$ the set of prime ideals 
  of $R$ of height $1$ and by $\cX=\cX_A$ the full subcategory of $A\Md$ consisting of all modules $M$ such that 
  $M_\gP=0$ for all $\gP\in\kP$. 
  We also consider two full subcategories of $A\Md$ related to $\kX$:
   \begin{align*}
    \kM=\kM(A)&=\{M\in A\Md\mid  X\in\kX\,\Arr\, \Hom_A(X,M)=0\} \\
    				&=\{M\in A\Md\mid X\in\kX\, \&\, X\sbe M \,\Arr\, X=0\}
      \end{align*}
      and
    {\small  \[    				
   \kC=\kC(A)=\{M\in A\Md\mid  X\in\kX\,\Arr\, \Hom_A(X,M)=\Ext^1_A(X,M)=0\}.
  \]} 
  Following P.\,Gabriel \cite{gab}, we call modules from $\kC$ \emph{closed}. Note that every injective module from $\kM$
  is closed. 
  Obviously, $\kX$ is a Serre subcategory of $A\Md$, so the quotient category $A\tMd$ is defined. Moreover, every $A$-module
  $M$ contains the biggest submodule from $\kX$, namely $M_{\kX}=\sum_{\substack{N\sbe M\\N\in \kX}}N$. Therefore, $\kX$ is a
  \emph{localizing subcategory} in $A\Md$, that is the natural functor $F:A\Md\to A\tMd$ has a right adjoint $G$  
  \cite[p.\,375,\,cor.\,1]{gab}. We denote  by $\thom_A(M,M')$ the set of morphisms $FM\to FM'$ in the category $A\tMd$.
  If $\al$ is such that $F\al$ is an isomorphism, we say, following \cite{bourbaki} that $\al$ is a \emph{pseudo-isomorphism}.
  It means that all localizations $\al_\gP\ (\gP\in\kP)$ are isomorphisms or, equivalently, that both $\ker\al$ and $\cok\al$ are in $\cX$.
  
  The following facts are those from \cite[Ch.\,III]{gab}.
  
   \begin{fact}\label{f1} 
  \begin{enumerate}
   \item  $FG\simeq \id_{A\tMd}$.
   \item If $M$ is closed, for arbitrary module $N$ the functor $F$ induces an isomorphism 
   $\Hom_A(N,M)\simeq\thom_A(N,M)$.
   \item $M$ is closed iff the adjunction map $M\to GFM$ is an isomorphism. Equivalently, $M\simeq GN$ for some $N\in A\tMd$.
   \item The functors $F$ and $G$ induce an equivalence of the categories $\kC(A)$ and $A\tMd$.
   \item  Both $\kX$ and $A\tMd$ (hence $\kC$) are \emph{Grothendieck categories}.
 \end{enumerate} 
   \end{fact}
   
    \begin{exam}\label{ex1} 
    Let $R$ be local noetherian of Krull dimension $2$. 
    \begin{itemize}
    \item $\kM(R)$ consists of modules that have no simple (hence no artinian) submodules.
    \item $\kC(R)$ is the subcategory of maximal Cohen-Macauley modules (maybe infinitely generated).
    \end{itemize}
    \end{exam}  
  
    On the other hand, $\kM$ is a reflective subcategory, that is the embedding functor $\kM\emb A\Md$ has a left adjoint,
  namely the functor $M\mapsto M/M_\kX$. Following I.\,Beck \cite{beck}, we call modules from $\kM$ \emph{codivisorial}.
  The following considerations explain this terminology.
  
   \begin{dfn}\label{defdiv} 
   Let $M$ be a torsion free $R$-module (for instance an ideal of $R$). We set
   \[
    M_\kP=\{x\in QM\mid \forall\,\gP\in\kP\ \exists\,r\in R\=\gP\ rx\in A\}
   \]
   and call $M$ \emph{divisorial} if $M_\kP=M$. We call the ring $R$ \emph{\pk} if it is divisorial as $R$-module.
   Krull rings are just \pk\ normal domains. A \emph{\pko} over $R$ is a \pno\ which is divisorial as $R$-module.
   
   \smallskip
   Note that if $A$ is not \pk\ and $\gA\sb A$ is a left (or right) ideal, it can happen that $\gA_\kP\not\sbe A$. We call $\gA$
   \emph{\qd} if $\gA=A\cap\gA_\kP$. If $A$ is \pk, \qd\ is the same as divisorial.
   \end{dfn}
     
     \begin{prop}\label{11} 
    $M\in\kM$ iff for every nonzero $x\in M$ the left ideal $\gA=\ann_Ax$ is \qd. 
     \end{prop}
      \begin{proof}
      If $a\in\gA_\kP$, for every $\gP\in\kP$ there is an element $r\in R\=\gP$ such that $ra\in\gA$, hence $r\in\ann_Rax$.
      Therefore $(Aax)_\gP=0$ for every $\gP\in\kP$. If $M\in\kM$ it implies that $ax=0$, that is $a\in \gA$ and $\gA$ is
      divisorial. On the other hand, if $M\notin\kM$, there is a nonzero element $x\in M$ such that $\frac{x}{1}=0$ in every 
      localization $M_\gP$, where $\gP\in\kP$, that is there is $r\in R\=\gP$ such that $rx=0$. For the ideal $\gA=\ann_Rx$ 
      it means that $\gA\cap R\not\sbe\gP$, whence $\gA_\gP=A_\gP$ and $\gA_\kP=A_\kP$, so $\gA$ is not \qd.
      \end{proof}
    
     \begin{cor}\label{12} 
     \begin{enumerate}
     \item  If $M$ is codivisorial, the ideal $\ann_AM$ is \qd.
     \item  An ideal $\gA\sbe A$ is \qd\ iff the module $A/\gA$ is codivisorial.
     \end{enumerate}
     \end{cor}
     
   \begin{prop}\label{13} 
    \begin{enumerate} 
    \item  If $M',M''$ are codivisorial and $0\to M'\to M\to M''$ is an exact sequence, $M$ is codivisorial.
    \item  If $M\emb M'$ is an essential extension of a codivisorial module $M$, then $M'$ is codivisorial. 
    In particular, the injective envelope of a codivisorial module is codivisorial (hence closed).
    \item  If $M\sbe M'$, $M'$ is codivisorial and $M'/M\in \kX$, then $M'$ is an essential extension of $M$.
    \end{enumerate}
   \end{prop}
    \begin{proof}
    (1) and (2) are obvious.
    
    (3) Let $x\in M'\=M$, $\ox=x+M\in M'/M$, $\gA=\ann_R x$ and $\gB=\ann_R\ox$. As $M'/M\in\kX$, the ideal $\gB$
    is contained in neither ideal $\gP\in\kP$. On the other hand, as $M'\in \kM$, there is a prime ideal $\gP$ such that
    $\gA\sbe \gP$. Hence $\gB\ne\gA$ and there is $a\in R$ such that $ax\ne0$ but $ax\in M$. It means that the
    extension $M\sbe M'$ is essential.
    \end{proof}
  
   We denote by $E(M)$ (or $E_A(M)$, if necessary) the injective envelope of the $A$-module $M$.
   
    \begin{cor}\label{14} 
       Let $M'=M/M_\kX,\ E=E(M')$ and $c(M)$ be the preimage in $E$ of $(E/M')_\kX$. Then $c(M)\simeq GFM$.
    \end{cor}
       \begin{proof}
      As $FM\simeq FM'$, we may suppose that $M$ is codivisorial. Then the adjunction map $\ga:M\to GFM$ is a 
      monomorphism and $\cok\ga\in\kX$. Therefore, $\ga$ is an essential monomorphism, thus there is a monomorphism
      $\ga':GFM\to E$ such that $\im\ga'/M\in\kX$, that is $\im\ga'\sbe c(M)$. Moreover, $c(M)/\im\ga\in\kX$. Since
      $\im\ga'\simeq GFM$ is closed, $\im\ga'=c(M)$.
      \end{proof}
          
      Recall the following well-known facts concerning rings and modules of fractions.
      
       \begin{fact}\label{f2} 
       Let $S$ be a multiplicative subset of $R$.
       \begin{enumerate}
       \item The natural embedding $\ms A\Md\to A\Md$ is fully faithful.
       \item  The natural map $\ms\Hom_A(M,N)\to\Hom_{\ms A}(\ms M,\ms N)$ is injective if $M$ is finitely generated
       and bijective if it is finitely presented.
       \item  An $\ms A$-module $M$ is injective iff it is injective as $A$-module.
       \item  If $E^*(M)$ is a minimal injective resolution of an $\ms A$-module $M$, it is also its minimal injective resolution
       as of $A$-module. 
       \end{enumerate}
       \end{fact}
     
      \begin{cor}\label{15} 
   If $M$ is an $A_\gP$-module or a $Q$-module and $X\in\kX$, then $\Ext^i_A(X,M)=0$ for all $i$. In particular, $M$ is closed.
      \end{cor}
      
       \begin{thm}\label{t1} 
        \begin{enumerate}
        \item If $M$ is torsion, $GFM\simeq\coprod_{\gP\in\kP} M_\gP$.
        \item If $M$ is torsion free, $GFM\simeq M_\kP$.
        \end{enumerate}
       \end{thm}
        \begin{proof}
        (1) As $M$ is torsion and $R$ is \pn, the image $M'$ of the natural map $\ka:M\to\prod_{\gP\in\kP}M_\gP$ belongs to
        $\coprod_{\gP\in\kP}M_\gP$. Morover, $\ka_\gP$ is an isomorphism for every $\gP\in\kP$. Therefore,
        $F\ka$ is an isomorphism as well  as $GF\ka$. As all $M_\gP$ are closed, it implies the claim.
        
        (2) As $Q$ is semisimple, $QM$ is injective, thus $QM=E(M)$. Note that $M_\gP=(M_\kP)_\gP$
        for every $\gP\in\kP$, whence $c(M)=c(M_\kP)$. By definition of $M_\kP$, it is the union of all submodules $N\sbe QM$
        such that $N\spe M$ and $(N/M)_\gP=0$ for all $\gP\in \kP$. Therefore, $M_\kP=c(M)\simeq GFM$.
        \end{proof}
          
     \begin{remk}\label{pk} 
     Note that every divisorial (in particular closed) $A$-module is actually an $A_\kP$-module. In particular,
     $A\tMd\simeq A_\kP\tMd$. Therefore, studying the category $A\tMd$ and divisorial modules, we may always
     suppose that $A$ is a \pko.
     \end{remk}
   
      \section{Injective modules and global dimension}
   \label{s2}     
   
   Now we describe codivisorial (or, the same, closed) injective modules or, equivalently, injective objects in the category $A\tMd$. 
   As mentioned in Rem.\,\ref{pk}, we may (and will) suppose that the ring $A$ is \pk. First, we note some general facts about prime 
   ideals of algebras. 
   
    \begin{fact}\label{p21} 
    Let $A$ be an $R$-algebra, $\gP$ be a prime ideal of $R$ such that $R_\gP$ is noetherian and $A_\gP$ is a finite $R_\gP$-algebra.
    Set $$\gP^\uparrow=\{\dP\mid \dP \text{\emph{ is a prime ideal of $A$ such that }} \dP\cap R=\gP\}.$$
    \begin{enumerate}
    \item  $\gP^\uparrow$ is nonempty and finite. 
    \\[.5ex]
    \emph{Let now $\dP\in\gP^\uparrow$ and $\gR_\gP=\rad A_\gP$.}
    \vskip.5ex
    \item  If $\dP'\in\gP^\uparrow$ and $\dP\sbe\dP'$, then $\dP=\dP'$.
    \item  $\dP_\gP$ is a maximal ideal in $A_\gP$ and every maximal ideal of $A_\gP$ is of this form.
    \item  $\gR_\gP=\bap_{\dP\in\gP^\uparrow}\dP_\gP$.
    \item  There is a unique simple $A_\gP$-module $U_\dP$ such that $\ann_AU_\dP=\dP$ and 
    $A_\gP/\dP_\gP\simeq U_\dP^{m(\dP)}$ for some $m(\dP)$.
    \item  $E(A/\dP)\simeq E(U_\dP)^{m(\dP)}$.
    \end{enumerate}
    \end{fact}
     \begin{proof}
     The proofs easily follow from \cite[Sec.\,5.6]{gab} and \cite[Sec.\,3.1]{morita} and the obvious remark that a prime ring
     has no zero divisors in its center.
     \end{proof}
    
   From now on $A$ is a reduced \pno\ over a reduced \pn\ ring $R$. We denote by $\kP(A)$ the set $\bup_{\gP\in\kP(R)}\gP^\uparrow$.
   
    \begin{thm}\label{it1} 
    Let $E$ be an indecomposable injective codivisorial $A$-module. Then either $E\simeq Q_j$ for some $j$ or $E\simeq E(U_\dP)$
    for some $\dP\in\kP(A)$.
    \end{thm} 
     \begin{proof}
  Note that $E= E(M)$, where $M=Ax$ for any nonzero element $x\in E$. It implies that $E$ is either torsion or torsion free. 
  If it is torsion free, then the natural map $M\to QM$ is an essential embedding. As $QM$ is injective, $E\simeq QM\simeq Q_j$ 
  for some $j$. Let $M$ be torsion, $\gA=\ann_Rx$. As $\gA$ contains \nzd{s}, the set $\kA=\{\gP\in\kP(R)\mid \gP\spe\gA\}$ is
  finite: $\kA=\{\lst\gP s\}$. We have seen in the proof of Thm.\,\ref{t1} that $\ka:M\to \bop_{i=1}^s M_{\gP_i}$ is an 
  essential embedding. As $E(M)$ is indecomposable, $E(M)\simeq E(M_{\gP_i})$ for some $i$.
  As $M_{\gP_i}$ is a torsion module over $A_\gP$, it contains a simple submodule, which is isomorphic to $U_\dP$ for
  some $\dP\in \kP(A)$. Therefore, $E\simeq E(U_\dP)$. 
   \end{proof}
   
   \begin{lem}\label{il1} 
     Let $f:\gA\to E(U_\dP)$, where $\gA$ is an ideal in $A$ and $\dP\in\gP^\uparrow$. If $f\ne0$ then $\ker f\cap R\sbe\gP$. 
      \end{lem}
       \begin{proof}
       Extend $f$ to a homomorphism $f':A\to E(U_\dP)$. Let $\gK=\ker f'$. Then $A/\gK$ embeds into $E(U_\dP)$,
       hence contains a submodule isomorphic to $U_\dP$. Therefore, $\gK\cap R\sbe\ann_AU_\dP\cap R=\gP$.
       \end{proof}
       
     \begin{thm}\label{it2} 
    Any coproduct of indecomposable injective codivisorial modules is injective and any codivisorial injective module is a
    coproduct of indecomposable modules.
    \end{thm}
     The proof consists of several claims.
      \begin{claim}
      Every multiple $E(U_\dP)^{(I)}$, where $\dP\in\kP(A)$, is injective.
      \end{claim}\label{it21} 
       \begin{proof}
       It is a coproduct of injective modules over the noetherian ring $A_\gP$, where $\gP=\dP\cap R$. 
       Hence it is injective \cite{matlis}.
       \end{proof}
       
        \begin{claim}
       Every coproduct of indecomposable injective codivisorial modules is injective. 
        \end{claim}\label{it22} 
         \begin{proof}
         It is true for every coproduct of modules $Q_j$, which is a module over the semisimple ring $Q$. So we have to
         prove it for a coproduct $P$ of modules $E(U_\dP)$, where $\dP$ runs through $\kP(A)$. Let $P(\gP)$ 
         be the part of this coproduct consisting of all $E(U_\dP)$ with $\dP\in\gP^\uparrow$ for a fixed $\gP$, so that
         $P=\coprod_{\gP\in\kP}P(\gP)$. As $\gP^\uparrow$ is finite and all sums of $E(U_\dP)$ are injective, each
         $P(\gP)$ is injective as well as every finite sum of them. Let $f:\gA\to P$, where $\gA$ is an ideal of $A$. As $P(\gP)$
         is closed and torsion and $\Hom_A(M,N)=0$ if $M$ is an $A_\gP$-module and $N$ is an $A_{\gP'}$-module with
         $\gP'\ne\gP$, $f$ factors as
         \[
          \gA \xarr{\,\pi} \gA/\ker f \xarr{\,\ka\,} \coprod_{\gP} \gA_\gP/\ker f_\gP \xarr{\, \vi\,}\coprod_{\gP}  P(\gP),
         \]
         where $\ka$ induces the isomorphism from Thm.\,\ref{t1}(1) and $\vi=\coprod_{\gP}\vi_\gP$. We also have a 
         commutative diagram
         \[ 
           \xymatrix@R=1em{ \gA \ar[r]^{\pi\ \ \ } \ar[d] & \gA/\ker f  \ar[r]^{\ka\ \ \ } \ar[d] & \coprod_{\gP} \gA_\gP/\ker f_\gP  \ar[d] \\
            A \ar[r] & A/\ker f  \ar[r] & \coprod_{\gP} A_\gP/\ker f_\gP            }
         \]
        where all vertical maps are embeddings. As $P(\gP)$ is injective, each map $\vi_\gP$ extends to a homomorphism
        $A_\gP/\ker f_\gP\to P(\gP)$. It gives an extension of $f$ to a homomorphism $A\to \coprod_{\gP}  P(\gP)$. 
        Therefore, this coproduct is injective.
        \end{proof}
        
          \begin{claim}\label{it23}
         Every codivisorial injective module $E$ is a coproduct of indecomposable modules.
          \end{claim}
           \begin{proof}
           Consider submodules of $E$ which are coproducts of indecomposables. Zorn lemma guarantees that there is a
           maximal $E'$ among them. Then $E=E'\+E''$ for some $E''$. If $0\ne x\in E''$, then $E''$ has a direct summand
           $E(Ax)$. As we have seen in the proof of Thm.\,\ref{it1}, $E(Ax)$ has a direct summand $I$ isomorphic to some $Q_j$
           or to some $E(U_\dP)$. As $E'\+I$ is bigger than $E'$, it is impossible. Hence $E'=E$.
		  \end{proof}
      
      \begin{cor}\label{ic2} 
    If all modules $M_i\ (i\in I)$ are codivisorial, $E(\coprod_{i\in I}M_i)\simeq\coprod_{i\in I}E(M_i)$
    \end{cor}

  For a divisorial torsion free $A$-module $M$ and an ideal $\gP\in\kP$ we denote 
  $$M_\gP^+=\{x\in (QM)_\gP\mid \gR_\gP x\sbe M_\gP\}\simeq\Hom_A(\gR_\gP,M_\gP).$$
  Obviously, $M^+_\gP/M_\gP=\soc_{A_\gP}(QM)_\gP/M_\gP\simeq \soc_A(QM/M)_\gP$ and
  $M^+_\gP/M_\gP\simeq \bop_{\dP\in\gP^\uparrow}U_\dP^{(r_M(\dP))}$ for some cardinalities
  $r_M(\gP)$.  
     
    \begin{cor}\label{ic2} 
    If $M$ is a torsion free divisorial $A$-module,
    $$ E(QM/M)\simeq \coprod_{\dP\in\kP(A)} E(U_\dP)^{(r_M(\dP))}. $$    
    \end{cor}    
  
     \begin{lem}\label{il2} 
         Let $M$ be codivisorial.
        \begin{enumerate}
        \item $\al:M\emb M'$ is an essential embedding  iff so is $T\al:TM\to TM'$.
        \item If $M$ is injective, so is $TM$.
        \item $E(TM)\simeq TE(M)$.
        \item Injective modules in $A\tMd$ are just coproducts of copies of $TQ_\gP$ and $TE(U_\dP)$, where $\dP\in\kP(A)$. 
        \end{enumerate}
     \end{lem}   
      \begin{proof}
      (1) is \cite[p.\,373,\,Lem.\,3]{gab} and (2) is \cite[p.\,374,\,Prop.\.6]{gab}.
      
      (3) follows immediately from (1) and (2).
      
      (4) follows from (3) and Thms.\,\ref{it1} and \ref{it2}, since every injective object is an injective envelope of each of
      its subobjects and every object in $A\tMd$ is isomorphic to $TM$ for some codivisorial $M$.
      \end{proof}  

     \begin{thm}\label{it3} 
      \begin{enumerate}
      \item $\idim TM=\sup_{\gP\in\kP}\idim_{A_\gP}M_\gP$.
      \item $\gdim A\tMd=\sup_{\gP\in\kP}\gdim A_\gP$.\!%
      \footnote{\,Note that left and right global dimensions of $A_\gP$ are equal since it is noetherian.}
      \end{enumerate}
     \end{thm}
      \begin{proof}
      (1) As $TM\simeq T(GFM)$, we can suppose that $M$ is closed. If $M$ is torsion, it is isomorphic to
      $\coprod_{\gP\in\kP}M_\gP$ and its minimal injective resolution is a coproduct of minimal injective resolutions
      of the modules $M_\gP$. Since all these modules are closed, it implies (1) by Lem.\,\ref{il2} and Fact\,\ref{f2}(4).
      
      Let now $M$ be torsion free. Then $E(M)=QM$. $M$ is injective iff $M=QM$, iff $M_\gP=Q_\gP$ for all
      $\gP\in\kP$. Otherwise, $\idim M=\idim QM/M+1$. As we suppose that $M$ is closed, 
      $TM$ is also not injective. Therefore, also $\idim TM=\idim T(QM/M)+1$. As $QM/M$ is torsion, we have already seen that
      $\idim T(QM/M)=\sup_{\gP\in\kP}\idim_{A_\gP}(QM/M)_\gP$. As $(QM/M)_\gP\simeq QM_\gP/M_\gP$ and
      $QM_\gP$ is the injective envelope of $M_\gP$, we have the formula (1) for torsion free $M$ too.
      
      If $M\in\kM$ is arbitrary, there is an exact sequence $0\to M'\to M\to M''\to 0$, where $M'$ is torsion and $M''$
      is torsion free. As $\idim M=\sup\{\idim M',\idim M''\}$ and (1) is valid for both $M'$ and $M''$, it is valid for $M$.
      
      (2) is an immediate corollary of (1).
      \end{proof}
      
       \begin{exam}\label{ex2} 
      Let $R$ be a local noetherian ring of Krull dimension $2$, $A$ be a finite $R$-algebra which is a maximal
      Cohen-Macaulay $R$-module. Then $\kC(A)=\mathrm{MCM}(A)$, the category of $A$-modules which are 
      maximal Cohen-Macaulay as $R$-modules (not necessarily finitely generated). Therefore, 
      $\mathrm{MCM}(A)\simeq A\tMd$ is a Grothendieck category and $\gdim\mathrm{MCM}(A)=\sup_\gP\gdim A_\gP$. 
      In particular, this category is \emph{hereditary} (of global dimension $1$)  iff all localizations $A_\gP$ are hereditary, 
      for instance, $R$ is normal and $A=R$ or $A$ is a subring of $\mathrm{Mat}(n,R)$ consisting of 
      matrices $(a_{ij})$ such that $a_{ij}\in\wrd\gP k$ for $i<j$, where $\lst\gP k$ are different prime ideals of height $1$.
       \end{exam}
 
     \begin{thm}\label{it4} 
          Let $M$ be a torsion $A$-module.
     \begin{enumerate}
     \item  If $E^*(\gP)$ is a minimal injective resolution of $M_\gP$, where $\gP\in\kP$, then $F\Si E^*$ is a minimal injective 
     resolution of $FM$, where $\Si E^*=\coprod_{\gP\in\kP}E^*(\gP)$.
     \item  If $N$ is a finitely generated $A$-module, then $$\tExt_A^i(N,M)\simeq \coprod_{\gP\in\kP}\Ext_A(N_\gP,M_\gP).$$
     \end{enumerate}
     \end{thm}
      \begin{proof}
       By Thm.\,\ref{t1}, $FM\simeq \coprod_{\gP\in\kP}FM_\gP$.  A minimal injective resolution of $M_\gP$ as of $A$-module
      coincides with its minimal injective resolution as of $A_\gP$-module, hence consists of torsion injective $A_\gP$-modules, 
      which are sums of some $E(U_\dP)$. By Lem.\,\ref{il2}, $FE^*(\gP)$ is a minimal injective resolution of $FM_\gP$. 
      As direct sums of injective codivisorial modules are injective, it implies (1).
      
       By Fact\,\ref{f1}(2), 
       $$\thom_A(FN,F\Si E^*)= \Hom_A(N,\Si E^*)= \coprod_{\gP\in\kP}\Hom_A(N,E^*(\gP))$$ 
       since $N$ is finitely  generated. Taking cohomologies, we obtain (2).
      \end{proof}
      
    We call a \pko\ $A$ \emph{pseudo-hereditary} if all localizations $A_\gP\ (\gP\in\kP)$ are hereditary. So are,   
    for instance, Krull rings. 
      
       \begin{cor}[Cf.\,{\cite[Ch.\,VII,\,\S\,1,\,Thm.\,4]{bourbaki}}]\label{ic4}
       Let $A$ be a pseudo-here\-di\-tary order and $M$ be an $A$-module.
        \begin{enumerate}
         \item  If $\tM$ is finitely generated, $FM\simeq F(\trs M)\+F(\tf M)$. 
         \item  If $M$ is finitely generated itself, there is a pseudo-isomorphism $f:M\to \trs M\+\tf M$.
         \end{enumerate} 
       \end{cor}
        \begin{proof}
        (1) follows from Thm.\,\ref{it4}, since $(\tf M)_\gP$ is a projective $A_\gP$-module.
        
        (2) There is a commutative diagram with exact rows
        \[
         \xymatrix{ 0\ar[r] & M' \ar[r]^u \ar[d]_\al & M \ar[r]^v \ar[d]_\be & M'' \ar[r] \ar[d]_\ga & 0 \\
                			  0\ar[r] & \tM' \ar[r]^{\ti u}  & \tM \ar[r]^{\ti v}  & \tM'' \ar[r]  & 0},
        \]      
        where $\al,\be,\ga$ are pseudo-isomorphisms.
        Here $M'=\trs M$, $M''=\tf M$ and we write $\ti X$ instead of $GFX$. As the exact sequence
        $0\to FM'\to FM\to FM'\to 0$ splits, the lower row of this diagram splits too, so there are morphisms $u':\tM\to \tM'$
        and $v':\tM''\to\tM$ such that $u'\ti u$ and $\ti vv'$ are identity maps. Then $\be=\ti uu'\be+v'\ga v$ and
        $\al=u'\be u$:
        \[
         \xymatrix{ 0\ar[r] & M' \ar[r]^u \ar[d]_\al & M \ar[r]^v \ar[d]_\be \ar[dl]_{u'\be} \ar[dr]^{\ga v}
         						& M'' \ar[r] \ar[d]^\ga & 0 \\
                			  0\ar[r] & \tM' \ar[r]^{\ti u}  & \tM \ar[r]^{\ti v} \ar@/^/[l]^{u'} & \tM'' \ar[r] \ar@/^/[l]^{v'}  & 0}
        \]
        Recall that $\tM'\simeq\prod_{\gP\in\kP}M'_\gP$. As $M$ is finitely generated, actually $u'\be$ maps $M$ to
        a finite sum $\bop_{i=1}^kM'_{\gP_i}$ for some $\lst \gP k$. One easily sees that 
        \[
         \hom_A(M,{\bop}_{i=1}^kM'_{\gP_i})\simeq \hom_A(M,S^{-1}M')\simeq S^{-1}\hom_A(M,M'),
        \]
        where $S=A\=\bup_{i=1}^k\gP_i$. It implies that $su'\be=\al\th$ for some $\th:M\to M'$ and some $s\in S$. 
        Note that multiplication by $s$ is an automorphism of $\tM'$, hence the morphism $f:M\to \tM'\+\tM''$ with the
        components $su'\be$ and $\ga v$ is a pseudo-isomorphism. The calculations above show that $f$ actually
	    factors as
	    \[
	     M \xarr{\smtr{\th\\ v}} M'\+M'' \xarr{\smtr{\al &0\\0& \ga}} \tM'\+\tM''  \xarr{\smtr{\ti u & v'}} \tM.
	    \]
	    The second and the third morphism here are pseudo-isomorphisms, hence so is $\smtr{\th\\ v}:M\to M'\+M''$.
        \end{proof}
        
  Let $A$ be a \pno, $\gP\in\kP$ and $\hA_\gP$ be the $\gP$-adic completion of $A$. Then $\hA_\gP$ is semiperfect
  and all finitely generated torsion $A_\gP$-modules are actually $\hA_\gP$-modules. If $A$ is pseudo-hereditary, the
  last \emph{Corollary} from \cite{hmod} implies that there is a unique $A_\gP$-module $\dL(\dP,l)$ which is of length
  $l$ with $\dL(\dP,l)/\gR_\gP\dL(\dP,l)\simeq U_\dP$ for any given $l$ and $\dP\in\gP^\uparrow$. 
  
      \begin{cor}[Cf.\,{\cite[Ch.\,VII,\,\S\,1,\,Thm.\,5]{bourbaki}}]\label{ic5} 
       Let $A$ be a pseudo-here\-di\-tary \pko.
       \begin{enumerate}
       \item  For every $\dP\in\kP(A)$ there is a unique up to isomorphism uniserial $A_\gP$-module $\dL(\dP,l)$ 
       of given finite length $l$ with the top $U_\dP$, where $\gP=\dP\cap R$. 
      \item  For every finitely generated torsion $A$-module $M$ there are prime ideals $\lst\dP k\in\kP(A)$, intergers
      $\lst rk$, $A$-modules $\lst Mk$ and pseudo-isomorphism $M\to\bop_{i=1}^kM_i$ such that  
      ${M_i}_\gP\simeq \dL(\dP,l)$, where $\gP=\dP\cap R$. The sequence of pairs $\{(\dP_i,l_i)\}$ is uniquely defined up
      to permutation.
       \end{enumerate}
        \end{cor}
         \begin{proof}
        (1) follows from the paper \cite{hmod} (the last Corollary).
          Then the proof of (2) is analogous to that of Cor.\,\ref{ic4}(2).     
         \end{proof}

\end{document}